\newif\ifplots
\plotstrue

\newif\ifelsevier
\elsevierfalse

\newif\ifieee
\ieeefalse

\newif\ifacc
\accfalse

\newif\ifflatdirs
\flatdirsfalse

\newif\ifdoublecolumn
\ifelsevier
    \doublecolumntrue
\else
    \ifieee
        \doublecolumntrue
    \else
        \ifacc
            \doublecolumntrue
        \else
            \doublecolumnfalse
        \fi
    \fi
\fi

\ifelsevier
    \documentclass[5p, times]{elsarticle}
    \def\pp{}
    \usepackage{natbib}
\else
    \ifieee
        \documentclass[letterpaper, 10 pt, conference]{ieeeconf}
        \def\pp{}
        \usepackage[sorting=none,citestyle=numeric,bibstyle=numeric,backend=bibtex,giveninits=true]{biblatex}
    \else
        \ifacc
            \documentclass[letterpaper, 10 pt, conference]{ieeeconf}
            \def\pp{}
            \usepackage[sorting=none,citestyle=numeric,bibstyle=numeric,backend=bibtex,giveninits=true]{biblatex}
        \else
            \documentclass[11pt]{article}
            \usepackage[margin=1.5in, headheight=42pt]{geometry}
            \def\pp{.}
            \usepackage[citestyle=alphabetic,bibstyle=alphabetic,backend=bibtex,giveninits=true]{biblatex}
        \fi
    \fi
\fi

\ifflatdirs
    \newcommand{\datadir}{.}
    \newcommand{\ext}{csv.txt}
\else
    \newcommand{\datadir}{data}
    \newcommand{\ext}{csv}
\fi

\usepackage{hyperref}

\usepackage{graphicx,psfrag,amsmath,versions}
\usepackage{array}
\usepackage{caption,subcaption}
\usepackage{tikz,pgfplots}
\usepackage{caption}
\usepgfplotslibrary{groupplots}
\usetikzlibrary{calc,patterns,decorations.pathmorphing,decorations.markings,arrows}
\usepackage{wasysym}
\usepackage[utf8]{inputenc}
\usepackage[english]{babel}
\usepackage{flushend}

\definecolor{green}{rgb}{0, .5, 0}
\definecolor{gray}{rgb}{.95, .95, .95}

\newcommand{\BEAS}{\begin{eqnarray*}}
\newcommand{\EEAS}{\end{eqnarray*}}
\newcommand{\BEQ}{\begin{equation}}
\newcommand{\EEQ}{\end{equation}}
\newcommand{\BIT}{\begin{itemize}}
\newcommand{\EIT}{\end{itemize}}

\newcommand{\eg}{{\it e.g.}}
\newcommand{\ie}{{\it i.e.}}

\newcommand{\reals}{{\mbox{\bf R}}}

\newcommand{\diag}{\mathop{\bf diag}}

\newcommand{\Expect}{\mathop{\bf E{}}}

\newcommand{\argmax}{\mathop{\rm argmax}}

\newcounter{oursection}

\ifelsevier
    \journal{Systems and Control Letters}
\fi

\title{Risk-Sensitive Model Predictive Control}

\author{Nicholas Moehle
}

\begin{document}

\ifelsevier
\else
    \maketitle
\fi

\begin{abstract}
We present a heuristic policy and performance bound for risk-sensitive convex stochastic control
that generalizes linear-exponential-quadratic regulator (LEQR) theory.
Our heuristic policy extends standard, risk-neutral model predictive control (MPC);
however, instead of ignoring uncertain noise terms,
our policy assumes these noise terms turns out either favorably or unfavorably,
depending on a risk aversion parameter.
In the risk-seeking case, this modified planning problem is convex.
In the risk-averse case, it requires minimizing a difference of convex functions,
which is done (approximately) using the convex--concave procedure.
In both cases, we obtain a lower bound on the optimal cost as a by-product of solving the planning problem.
We give a numerical example of controlling a battery to power an uncertain load,
and show that our policy reduces the risk of a very bad outcome
(as compared with standard certainty equivalent control)
with negligible impact on the the average performance.
\end{abstract}


\ifelsevier
    \maketitle
\fi

\section{Introduction}
In this paper, we study the problem of controlling a linear dynamical system driven by additive noise
in order to minimize a sum of convex stage costs, while satisfying state and control constraints.
In the standard \emph{risk-neutral} problem,
we minimize the expected value of this sum.
We focus on the \emph{risk-sensitive} problem,
in which we minimize the expected value of an exponential function of the cost.
This formulation is parameterized by a risk-aversion parameter $\gamma$.
For $\gamma > 0$, the problem is \emph{risk averse} or \emph{pessimistic};
for $\gamma < 0$, the problem is \emph{risk seeking} or \emph{optimistic}.
This problem formulation is a generalization of the LEQR problem,
in which the stage costs are quadratic and the noise is Gaussian.

We give lower bounds on the optimal value of this problem
that are based on ideas from large deviations theory.
These bounds generalize the certainty equivalent bound (\ie, Jensen's inequality) 
obtained by solving an optimal planning problem that replaces the additive noise term with its expected value.

Evaluating our bound requires solving an optimization problem,
which we use as the basis for a control policy
we call \emph{risk-sensitive model predictive control} (RS-MPC).
As opposed to other LEQR extensions in literature,
RS-MPC handles non-smooth convex stage cost functions
(which can encode convex state and control constraints)
as well as non-Gaussian disturbances.
In the risk-averse case, evaluating the RS-MPC policy
requires solving a minimax problem in which we plan against a worst-case disturbance;
in the risk-seeking case, we co-optimize over the disturbance 
along with the control and state trajectories.

\subsection{Related work}

\paragraph{Certainty equivalence for LEQR\pp}
The basic \emph{linear quadratic regulator} (LQR) problem is 
to control a linear dynamical system with an additive disturbance to minimize the expected value of a sum of quadratic stage costs.
The \emph{certainty equivalence principle} (CEP)
states that ignoring the stochastic noise, solving the optimal planning problem,
and then applying the optimal first input results in an optimal control policy \cite[\S 3.1]{bertsekas2017dynamic}.
Furthermore, the planned state and input trajectories describe the mean trajectories under such an optimal policy.

The LEQR problem swaps out the expectation operator for a risk-sensitive certainty equivalent operator,
\ie, we minimize the expected value of an exponential function of the total cost.
Whittle describes a \emph{risk-sensitive} certainty equivalence principle (RS-CEP) for LEQR,
in which the deterministic planning problem is a two-player game
between the planner and ``nature'' \cite[\S 10.2]{whittle1990risk}.
For risk-\emph{averse} LEQR, this game is adversarial,
while for risk-\emph{seeking} LEQR, it is cooperative.
More specifically, nature chooses a value of the disturbance that trades off pessimism (or optimism) with plausibility,
and the planner optimizes accordingly.
The (risk-neutral) CEP for LQR is the special case in which we are not optimistic or pessimistic,
and therefore nature selects the most plausible values for the disturbance.
(An example of a similar risk-averse CEP can be found in \cite{moehle2021certainty}.)

\paragraph{Model Predictive Control\pp}
Model predictive control (MPC)
is a heuristic technique that applies the certainty-equivalence principle
beyond where it is theoretically justified,
\eg, to problems with non-quadratic stage cost functions
\cite[\S 4.3]{bertsekas2017dynamic}.
An MPC policy replaces all uncertain quantities with estimates,
then solves the resulting (deterministic) optimal planning problem.
This is not optimal in general, but typically yields excellent practical performance.
In some contexts, MPC is also called
\emph{certainty-equivalent control} or \emph{receding-horizon control};
see \cite{kwon2006receding, borrelli2017predictive}.

The method we propose in this paper (RS-MPC)
is similar in spirit to standard, risk-neutral MPC
in that it applies a CEP beyond where it is theoretically justified.
In our case, however, we apply the RS-CEP of LEQR
instead of the standard, risk-neutral CEP of LQR;
the resulting planning problem is a two-player game.
The RS-CEP policy can be fielded in much the same way as a (risk-neutral) MPC policy.

\paragraph{Iterative LEQR\pp}
Iterative LEQR is a heuristic for risk-sensitive nonlinear optimal control problems
that solves successive, local LEQR approximations of the problem around a candidate trajectory
\cite{farshidian2015risk, roulet2020convergence}.
A critical limitation of this approach
is the assumption that the stage cost functions are second-differentiable
and the state and control variables are unconstrained.
Our approach, while limited to linear dynamics, allows for non-smooth convex stage costs,
which can encode convex state and control constraints, as well as non-Gaussian disturbances.
Our focus on convexity also allows us to provide a global performance bound and convergence guarantee,
which are not possible using iterative LEQR.

\paragraph{Risk aversion and adversarial measures\pp}
Many results exist that equate risk-averse decision problem
with a zero-sum games in which an adversary chooses the probability measure
that the decision maker optimizes against.
(The most relevant for our case is \cite{petersen2000minimax}.)
In our approach, the adversary selects a \emph{specific value} of the disturbance,
which is typically a much more tractable problem than choosing a distribution.
(The cost of this tractability is that our game is not equivalent
but merely provides a bound on it.)


%


\ifieee 
\else
\subsection{Outline}
In section~\ref{s-risk},
we define our measure of risk, and we give an optimization-based bound on it.
In section~\ref{s-control-prob}, we define the risk-averse linear convex control problem.
We discuss the prescient relaxation of this problem in section~\ref{s-prescient},
and we use this relaxation as the basis for a heuristic policy.
In section~\ref{s-convex-concave}, we discuss the algorithmic details of the heuristic
policy in the risk-averse case.
We conclude with a numerical example in section~\ref{s-examples}.
\fi

\section{Risk}
\label{s-risk}
The \emph{risk} of a real-valued random variable $z$ is defined as
\begin{equation}
\label{e-risk-operator}
R_{\gamma}(z) = \frac1\gamma \log\Expect\exp(\gamma z),
\end{equation}
where $\gamma$ is the risk aversion parameter.
In this paper, $z$ represents a cost to be minimized,
and we call the case $\gamma > 0$ the \emph{risk-averse} case,
because it more heavily weights large values of $z$ than small values.
Likewise, the case $\gamma<0$ is \emph{risk seeking}.
We define $R_0(z) = \Expect z$;
we call this case \emph{risk neutral}.


\subsection{Risk bound}

\paragraph{Rate function\pp}
\label{s-rate-function}
The cumulant generating function $c:\reals^n\to\reals$ of a random vector $w$ is
\begin{equation}
\label{e-cgf}
c(y) = R_1(w^T y) = \log \Expect \exp w^T y.
\end{equation}
The cumulant generating function is convex, regardless of the distribution of $w$ \cite[pg. 106]{cvxbook}.
The \emph{rate function} $\rho:\reals^n\to\reals$
is the Fenchel conjugate of the cumulant generating function:
\[
\rho(x) = c^*(x) = \sup_y \big( x^T y - c(y) \big).
\]
The rate function appears in large deviations theory,
where it is used to approximate the distribution of the average of a large number
of independent samples of $w$. 
(See \cite[\S 18]{whittle2012probability} or \cite{den2008large} for an introduction.
Note that here we refer to the specific rate function defined in Cram\'{e}r's theorem,
as opposed to other rate functions that arise large deviations theory.)
The rate function can be interpreted
as a smoothed version of the negative log-likelihood function $\ell(x) = -\log p(x)$.
\ifieee
\else
In figure~\ref{f-rate-functions},
we compare the negative log-likelihood $\ell$ and the rate function $\rho$
for several common distributions.
\fi

It is easy to show that $\Expect w$ is the unique minimizer of $\rho$, and $\rho(\Expect w) = 0$.
(These properties derive from well-known properties of the cumulant generating function $c$,
as well as basic facts of convex analysis.)
Note that the cumulant generating function is the conjugate of the rate function,
\ie, $\rho^{*} = c$.

\ifieee
\else
    \begin{figure}
    \centering
    \ifplots
    \begin{tikzpicture}
        \begin{groupplot}[
            group style = {group size=2 by 2, vertical sep=1cm},
            table/col sep=comma,
            ymin = -0.5,
            ymax = 3,
            width = 0.5\columnwidth,
            height = 0.4\columnwidth,
            title style={yshift=-1.0ex},
            xticklabels={,,},
            yticklabels={,,},
        ]

        \nextgroupplot[
            xmin = -2,
            xmax = 2,
            title={Uniform},
        ]
            \addplot [blue] table{\datadir/ell_uniform.\ext};
            \addplot [green] table{\datadir/rho_uniform.\ext};

        \nextgroupplot[
            xmin = -2.5,
            xmax = 2.5,
            title={Gaussian},
        ]
            \addplot [blue] table{\datadir/ell_gaussian.\ext};
            \addplot [green] table{\datadir/rho_gaussian.\ext};

        \nextgroupplot[
            xmin = -7,
            xmax = 7,
            ymax = 7,
            title=Laplace,
        ]
            \addplot [blue] table{\datadir/ell_laplace.\ext};
            \addplot [green] table{\datadir/rho_laplace.\ext};

        \nextgroupplot[
            xmin = 0,
            xmax = 10,
            ymax = 5.5,
            title={${\rm Poisson}$}
        ]
            \addplot [blue, only marks, mark size=1pt] table{\datadir/ell_poisson.\ext};
            \addplot [green] table{\datadir/rho_poisson.\ext};

        \end{groupplot}
    \end{tikzpicture}
    \fi
    \caption{
    Rate functions $\rho(x)$ (in green) and (shifted) negative log-likelihood functions
    $\ell(x) - \ell(\Expect x)$ (in blue)
    for a uniform distribution,
    Gaussian distribution, 
    Laplace distribution 
    and Poisson distribution (with arrival rate $3$).
    }
    \label{f-rate-functions}
    \end{figure}
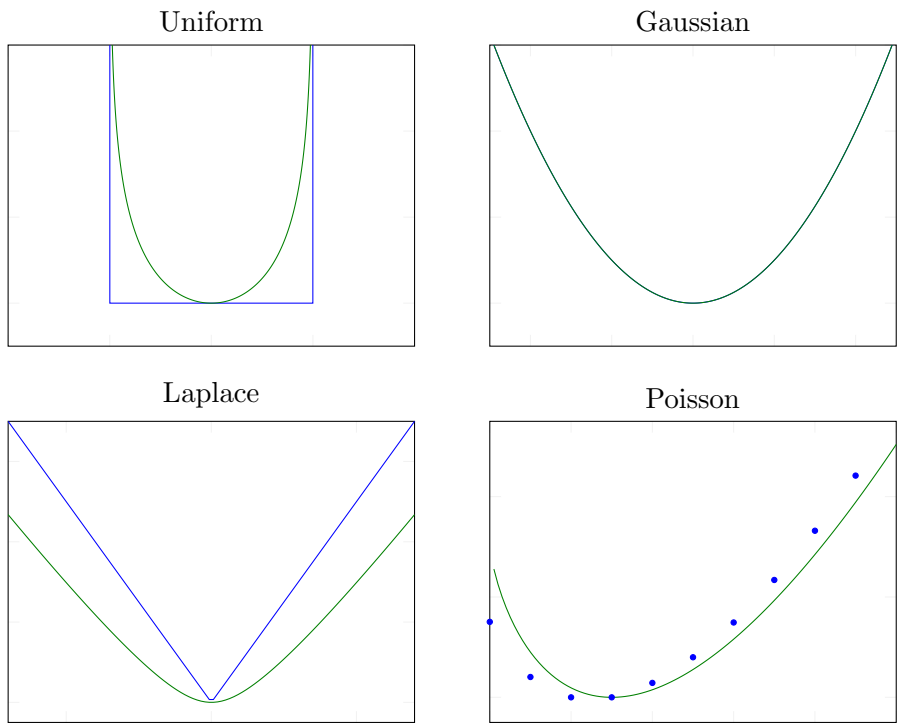
\fi


\paragraph{Risk bound\pp}
Consider a convex function $f:\reals^n\to\reals$ and a random variable $w\in\reals^n$ with rate function $\rho$.
For $\gamma \neq 0$, the following inequality holds:
\begin{align}
\label{e-varadhan-bound}
\frac1\gamma \sup_z \big( \gamma f(z) - \rho(z)\big) \le
R_\gamma \big(f(w)\big).
\end{align}
\ifieee
    This inequality is proven in the appendix.
\else
    This inequality is proven in appendix~\ref{s-varadhan-bound}.
\fi
It says that the value of $f$ at a single, well-chosen point $z$
approximates $R_\gamma(f(w))$
once adjusted for the likelihood of $z$, as measured by $\rho(z)$.

When $\gamma < 0$,
the quantity in the supremum is concave,
and evaluating the bound involves solving a simple convex optimization problem.
When $\gamma > 0$, we must instead maximize a difference of convex functions,
which is computationally hard in general;
we return to this issue in section~\ref{s-convex-concave}.

\paragraph{Comparison with Jensen's inequality\pp}
Take $\gamma > 0$.
Because $z=\Expect w$ is a valid choice in the left-hand side of \eqref{e-varadhan-bound},
and because $\rho(\Expect w) = 0$, we have
\begin{align*}
f(\Expect w) \le
\frac1\gamma \sup_z \big( \gamma f(z) - \rho(z)\big) \le
R_\gamma \big(f(w)\big),
\end{align*}
\ie, the bound given above is stronger than Jensen's inequality.
In fact, the bound~\eqref{e-varadhan-bound} reduces to Jensen's inequality in the risk-neutral case $\gamma\to 0$.
(This is because, in this limit, the the choice of $z$ in the supremum is dominated by $\rho$,
and because $z = \Expect w$ minimizes $\rho$ with the value $\rho(\Expect w) = 0$.)



\section{Risk-sensitive control}
\label{s-control-prob}

\paragraph{Dynamics\pp}
Consider the affine stochastic dynamics
\begin{equation}
\label{e-dynamics}
x_{t+1} = A_t x_t + B_t u_t + w_t, \quad t = 0, \ldots, T-1,
\end{equation}
defined over $T$ time periods.
Here $x_t\in\reals^n$ is the state, which has initial condition $x_0 = x_{\rm init}$,
and $u_t\in\reals^m$ is the control input.
The matrices $A_t$ and $B_t$ are deterministic.
The vectors $w_t\in\reals^n$ are random and independent across time periods,
with distributions $p_t$,
cumulant generating functions $c_t$, and rate functions $\rho_t$.

We use the compact notation 
\ifdoublecolumn
    \begin{align*}
    x &= (x_0, \dots, x_T)\in\reals^{n(T+1)} \\
    u &= (u_0, \dots, u_{T-1})\in\reals^{nT} \\
    w &= (w_0, \dots, w_{T-1})\in\reals^{nT}.
    \end{align*}
\else
    \begin{align*}
    x = (x_0, \dots, x_T), \quad
    u = (u_0, \dots, u_{T-1}), \quad
    w = (w_0, \dots, w_{T-1}),
    \end{align*}
\fi
and denote by $p$, $c$, and $\rho$
the probability distribution, cumulant generating function, and rate function of $w$.


\paragraph{Policy\pp}
A \emph{policy} $\pi$ is a function that maps the time period and state
to a control input, \ie, $u_t = \pi_t (x_t)$.

\paragraph{Cost\pp}
The total cost is defined as
\[
C_\pi(w) = g_T(x_T) + \sum_{t=0}^{T-1} g_t(x_t, u_t),
\]
where the stage cost functions $g_t$ are convex for all $t$.
We allow $g_t$ to be take the value $+\infty$,
which can be used to encode convex state and control constraints.
We emphasize that the total cost is a function of the policy $\pi$
as well as the random disturbance $w$ that obtains.
(The total cost is therefore a scalar-valued random variable).

\paragraph{Problem\pp}
The \emph{risk-sensitive linear convex control problem} is to choose a policy $\pi$ that minimizes
the risk-adjusted total cost:
\begin{equation}
\label{e-stoch-prob}
\begin{array}{ll}
\underset{\pi}{\mbox{minimize}} & J_\pi = R_\gamma\big(C_{\pi}(w) \big).
\end{array}
\end{equation}
We denote the infimum of $J_{\pi}$ over all policies as $J^\star$.

\ifieee
\else
    \paragraph{Breakdown\pp}
    The stochastic control problem may be unbounded ($J^\star = -\infty$)
    or infeasible ($J^\star  = \infty$).
    It may also be finite for some value of $\gamma$,
    but infinite for some larger value of $\gamma$.
    This phenomenon is called \emph{neurotic breakdown},
    and is simply a special interpretation of infeasibility
    due to a large value of $\gamma$.
    This may occur even if the stage cost functions $g_t$ only take finite values.
    It may also be that $J^\star$ is finite for some value of $\gamma$,
    but is $-\infty$ for some smaller value of $\gamma$.
    This is likewise called \emph{euphoric breakdown}.
\fi

\section{Prescient problem}
\label{s-prescient}
If the noise $w$ is known in advance,
the stochastic control problem reduces to the deterministic \emph{prescient problem}
\ifdoublecolumn
    \begin{equation} \label{e-prec-prob}
    \begin{array}{ll}
    \mbox{minimize} 
        & g_T(x_T) + \sum_{t=0}^{T-1} g_t(x_t, u_t) \\
    \mbox{subject to} 
        & x_{t+1} = A_tx_t + B_tu_t + w_t \\
        & x_0 = x_{\rm init}.
    \end{array}
    \end{equation}
    The variables are $x$ and $u$, and the first constraint holds for $t=0, \dots, T-1$.
\else
    \begin{equation} \label{e-prec-prob}
    \begin{array}{ll}
    \mbox{minimize} 
        & g_T(x_T) + \sum_{t=0}^{T-1} g_t(x_t, u_t) \\
    \mbox{subject to} 
        & x_{t+1} = A_tx_t + B_tu_t + w_t, \quad t=0, \dots, T-1 \\
        & x_0 = x_{\rm init}.
    \end{array}
    \end{equation}
    The variables are $x$ and $u$.
\fi
We denote by $C_{\rm pr}(w)$ the optimal value of \eqref{e-prec-prob} as a function of $w$,
and note that this function is convex.
In addition, $C_{\rm pr}(w)$ is random, because it depends on the random vector $w$.

\paragraph{Prescient bound\pp}
For any outcome $w$, prescient control obtains the lowest possible cost, \ie,
\[
C_{\rm pr}(w) \le C_{\pi}(w)
\]
for any policy $\pi$.
Because the risk operator $R_\gamma$ is monotonic,
we can apply it to both sides to obtain
\[
R_\gamma \big( C_{\rm pr}(w) \big) \le J_{\pi}.
\]
(Recall that the risk operator ``averages out'' the random variable $w$.)
By taking the infimum of the right-hand side over $\pi$, we obtain
\begin{equation}
\label{e-prec-bound}
R_{\gamma}\big(C_{\rm pr}(w) \big) \le J^\star,
\end{equation}
which says that the risk-adjusted value of the prescient problem is less than the optimal value of \eqref{e-stoch-prob}.
This bound is stronger than Jensen's inequality $C_{\rm pr}(\Expect w) \le J^\star$,
obtained by solving \eqref{e-prec-prob} with $w$ replaced its mean $\Expect w$.


\subsection{Bounds via rate function}
We now combine the risk bound~\eqref{e-varadhan-bound} with the prescient bound~\eqref{e-prec-bound},
taking $f = C_{\rm pr}$.
We do this separately for the risk-seeking case and the risk-averse case.

\paragraph{Risk-seeking case---Co-optimization over noise\pp}
For $\gamma < 0$, applying \eqref{e-varadhan-bound} to \eqref{e-prec-bound} and simplifying yields
\begin{equation}
\label{e-bound-risk-seeking}
\inf_w \Big( C_{\rm pr}(w) -\frac1\gamma \rho(w) \Big) \le J^\star.
\end{equation}
The left-hand side can be evaluated by solving the convex optimization problem
\ifdoublecolumn
    \begin{equation}
    \label{e-prec-prob-risk-seeking}
    \begin{array}{ll}
    \mbox{minimize} 
        & g_T(x_T) + \sum_{t=0}^{T-1} g_t(x_t, u_t) - (1/\gamma) \rho_t(w_t) \\
    \mbox{subject to} 
        & x_{t+1} = A_tx_t + B_tu_t + w_t \\
        & x_0 = x_{\rm init}
    \end{array}
    \end{equation}
    with variables are $x$, $u$, and $w$.
    The first constraint holds for $t=0, \dots, T-1$.
\else
    \begin{equation}
    \label{e-prec-prob-risk-seeking}
    \begin{array}{ll}
    \mbox{minimize} 
        & g_T(x_T) + \sum_{t=0}^{T-1} g_t(x_t, u_t) - (1/\gamma) \rho_t(w_t) \\
    \mbox{subject to} 
        & x_{t+1} = A_tx_t + B_tu_t + w_t, \quad t=0, \dots, T-1 \\
        & x_0 = x_{\rm init}
    \end{array}
    \end{equation}
    with variables are $x$, $u$, and $w$.
\fi
The optimal $w$ achieves the infimum in \eqref{e-bound-risk-seeking}
and the corresponding $x$ and $u$ are optimal for problem~\eqref{e-prec-prob}
with this value of $w$.

Problem~\eqref{e-prec-prob-risk-seeking} has the following interpretation.
In the risk-seeking case, we exhibit optimism,
\ie, we assume that the uncertain quantity $w$
will turn out in our favor.
In the resulting planning problem,
we co-optimize over the input, state, and noise trajectories.
We also ensure that $w$ is reasonably likely
by penalizing large values of $\rho(w)$.
(A similar phenomenon appears in the LEQR case; see \cite[\S 6.4]{whittle1990risk}.)

\paragraph{Risk-averse case---Adversarial noise\pp}
In the risk-averse case $\gamma > 0$,
applying \eqref{e-varadhan-bound} to \eqref{e-prec-bound} and simplifying yields
\begin{equation}
\label{e-bound-risk-averse}
C_{\rm pr}(w) - \frac1\gamma \rho(w) \le J^\star,
\end{equation}
which holds for any value of $w$.
This says that the value of problem~\eqref{e-prec-prob},
when adjusted to account for the likelihood of $w$,
is a lower bound on the optimal value of \eqref{e-stoch-prob}.

The tightest bound is obtained by maximizing the left-hand side over $w$.
In theory, this task is computationally difficult, as it involves maximizing over
the difference of two convex functions.
However, a very good heuristic, called the \emph{convex--concave procedure}, can be applied here,
and is discussed further in section~\ref{s-convex-concave}.
Furthermore, even suboptimal values of $w$ obtained by such a heuristic still produce a valid bound.

\paragraph{Risk-neutral case---Ignoring noise\pp}
As discussed in section~\ref{s-risk},
the bound~\eqref{e-bound-risk-averse} reduces to Jensen's inequality
as $\gamma \to 0$.
In the context of the linear-convex control problem, this results in the standard certainty equivalent bound
\[
C_{\rm pr}(\Expect w) \le J^\star.
\]

\subsection{Risk-sensitive certainty equivalent control}
\label{s-shrinking-horizon-policy}
Here we present RS-MPC, a heuristic policy based on the prescient problem~\eqref{e-prec-prob}.
To do this, we define how the control input $u_t$ is computed
as a function of the current state $x_t$ and time period $t$.
We will explain how to do this when $t = 0$ below.
To define the policy for $t = 1,\dots, T-1$,
we simply define a new stochastic control problem with initial state $x_t$
and horizon length $T - t$,
and then calculate the optimal first control input for this problem.
This approach is called \emph{shrinking-horizon control}
and is discussed in detail in \cite[\S 4.2]{skaf2010shrinking}.


\paragraph{Policy definition\pp}
We now define the initial input $u_0 = \pi_0(x_{\rm init})$.
In the risk-seeking case, we simply solve problem~\eqref{e-prec-prob-risk-seeking},
and use the optimal first control input $u_0$.
In the risk-averse case,
we carry out the following steps:
\begin{enumerate}
\item Find a maximizer $w^\star$ of $C_{\rm pr}(w) - (1/\gamma) \rho(w)$.
\item Solve \eqref{e-prec-prob} using $w = w^\star$,
and take $\pi_0(x_{\rm init})$ to be an optimal value of the first control input $u_0$.
\end{enumerate}
In the risk-averse case, this policy cannot be implemented exactly in practice,
because step 1 involves maximizing over a difference of convex functions,
which is a computationally hard problem.
The maximization in step 1 can instead be carried out approximately using the convex--concave procedure,
which is detailed in the next section.

\section{Convex--concave procedure}
\label{s-convex-concave}
We propose using the convex--concave procedure to find the best bound in \eqref{e-bound-risk-averse},
\ie, to approximately solve the problem
\begin{equation}
\label{e-max-over-w-prob}
\begin{array}{ll}
\mbox{maximize} 
    & \displaystyle C_{\rm pr}(w) - \frac1\gamma \rho(w)
\end{array}
\end{equation}
over the variable $w\in\reals^{nT}$.
In the RS-MPC policy of section~\ref{s-shrinking-horizon-policy},
this approximate method can be used in step 1 instead of carrying out the exact minimization over $w$.
For more information on the convex--concave procedure, see \cite{lipp2016variations}.

\subsection{Algorithm overview}
Starting with the initial guess $w^{(0)} = \Expect w$,
we define $w^{(k)}$ from $w^{(k-1)}$ by repeating the following steps.
\begin{enumerate}
    \ifdoublecolumn
        \item \emph{Minorization.} Form a first-order approximation 
            $\hat C_{\rm pr}(w; w^{(k-1)})$ of $C_{\rm pr}$ around $w^{(k-1)}$.
        \item \emph{Maximization.} Take
            \[
            \displaystyle w^{(k)} = \argmax_w \Big(\hat C_{\rm pr}\big(w; w^{(k-1)}\big) - \frac1\gamma \rho(w)\Big).
            \]
    \else
        \item \emph{Minorization.} Form a first-order approximation $\hat C_{\rm pr}(w; w^{(k-1)})$
        of $C_{\rm pr}$ around $w^{(k-1)}$.
        \item \emph{Maximization.} Take 
            $\displaystyle w^{(k)} = \argmax_w \Big(\hat C_{\rm pr}\big(w; w^{(k-1)}\big) - \frac1\gamma \rho(w)\Big)$.
    \fi
\end{enumerate}
We note that the objective of \eqref{e-max-over-w-prob},
evaluated at the iterates $w^{(k)}$, for $k = 0, 1, \dots$,
forms an increasing, convergent sequence \cite[\S 1.3]{lipp2016variations},
and can be used as a basis for a termination criterion.

%
%

\subsection{Implementation}
\label{s-implementation}
We now discuss implementation details of the algorithm,
which greatly simplify the algorithm steps.

\paragraph{Minorization step\pp}
To form a first-order approximation of $C_{\rm pr}$,
we require a subgradient of $C_{\rm pr}$ with respect to $w_t$, for $t=0,\dots,T-1$.
One such subgradient is an optimal dual variable $\lambda_t$
for the time-$t$ dynamics constraint of problem~\eqref{e-prec-prob}.
This means that a subgradient of 
$C_{\rm pr}(w)$ is $\lambda = (\lambda_0,\dots, \lambda_{t-1})\in\reals^{nT}$
and a first-order approximation of $C_{\rm pr}$ around $w'$ is
\[
\hat C_{\rm pr}(w, w') = C_{\rm pr}(w') + {\lambda}^T (w - w').
\]
Computing $C_{\rm pr}(w')$ and $\lambda$ requires solving problem~\eqref{e-prec-prob}.

\paragraph{Maximization step\pp}
The iterate $w^{(k)}$ maximizes
\[
C_{\rm pr}(w') + {\lambda}^T (w - w') - \frac1\gamma \rho(w)
\]
over $w$.
We drop the constant term $C_{\rm pr}(w') - {\lambda}^T w'$,
and instead maximize over ${\lambda}^T w - (1/\gamma)\rho(w)$.
The unique maximizing value of $w$ can be expressed in terms of the Fenchel conjugate of $\rho$,
which is the cumulant generating function $c$.
This maximizing value $w^\star$ is 
$w^\star = \nabla c(\gamma \lambda)$,
where $\nabla c$ is the gradient of the cumulant generating function of random variable $w$.

\subsection{Final, simplified algorithm}
\label{s-final-algo}
Starting with $w^{(0)} = \Expect w$, the iterates are defined as
\begin{enumerate}
\item \emph{Minorization.} Compute $\lambda^{(k-1)}$, the vector of optimal dual variables for problem~\eqref{e-prec-prob}
    with $w = w^{(k-1)}$.
\item \emph{Maximization.} Compute $\displaystyle w^{(k)} = \nabla c(\gamma \lambda^{(k-1)})$.
\end{enumerate}
We terminate the algorithm if the objective of \eqref{e-max-over-w-prob}, evaluated at $w^{(k)}$,
does not increase more than some positive value $\epsilon$
for a specified number of iterations.

\ifieee
\else

\section{LEQR}
As our first example, 
we revisit the classical linear-exponential-quadratic regulator problem.
In this case, we have
$g_T(x) = x^T Q x$ and 
\[
g_t(x, u) = x^T Q x + u^T R u, \quad t=0, \dots, T-1.
\]
where the matrices $Q$ and $R$ are positive semidefinite.
We also have $w_t \sim \mathcal N(0, \Sigma)$ for $t = 0, \dots, T-1$,
which means the rate function is 
\[
\rho(w) = \frac12 \sum_{t=0}^{T-1} w_t^T \Sigma^{-1} w_t.
\]
The prescient problem~\eqref{e-prec-prob} is a deterministic linear-quadratic control problem:
\begin{equation} \label{e-leqr-prob}
\begin{array}{ll}
\mbox{minimize} 
    & x_T Q x_T + \sum_{t=0}^{T-1} x_t^T Q x_t + u_t^T R u_t \\
\mbox{subject to} 
    & x_{t+1} = Ax_t + Bu_t + w_t, \quad t=0, \dots, T-1 \\
    & x_0 = x_{\rm init}.
\end{array}
\end{equation}
The optimal value $C_{\rm pr}(w)$ is a convex quadratic function of $w$.
As a result, the left-hand side of the bound~\eqref{e-varadhan-bound},
which is $C_{\rm pr}(w) - (1/\gamma) \rho(w)$, is also a quadratic function of $w$.
The maximizing value of $w$ can therefore be computed exactly, even in the risk-averse case.
(Indeed, in the risk-averse case,
the maximum value is finite if and only if this quadratic function is concave.)
Furthermore, the RS-MPC policy of section~\ref{s-shrinking-horizon-policy} is in fact optimal for LEQR.
This is discussed in \cite[\S 10]{whittle1990risk}.

In fact, it can be shown that this value of $w$,
as well as the corresponding optimal $x$ and $u$ for \eqref{e-prec-prob},
solve the system of linear equations
\[
\begin{bmatrix}
    \bar Q & 0 & \bar A^T & E_0 \\
    0 & \bar R & \bar B^T & 0 \\
    \bar A & \bar B & \gamma \bar \Sigma & 0\\
    E_0 & 0 & 0 & 0
\end{bmatrix}
\begin{bmatrix} x \\ u \\ w \\ \nu \end{bmatrix}
=
\begin{bmatrix} 0 \\ 0 \\ 0 \\ x_{\rm init} \end{bmatrix}
\]
\ifdoublecolumn
    where 
    \begin{equation*}
    \begin{split}
        \bar Q &= \diag(Q, \dots, Q), \\
        \bar B &= \diag(B, \dots, B),
    \end{split}
    \qquad
    \begin{split}
        \bar R &= \diag(R, \dots, R), \\
        \bar \Sigma &= \diag(\Sigma, \dots, \Sigma),
    \end{split}
    \end{equation*}
\else
    where $\bar Q = \diag(Q, \dots, Q)$, $\bar R = \diag(R, \dots, R)$,
    $\bar B = \diag(B, \dots, B)$, $\bar \Sigma = \diag(\Sigma, \dots, \Sigma)$,
\fi
and
\begin{align*}
\bar A = \begin{bmatrix} 
     A     &     -I &        &       &        \\ 
           &        & \ddots &       &        \\ 
           &        &        &     A &     -I \\ 
\end{bmatrix},
\quad
E_0 = \begin{bmatrix} I & 0 & \cdots & 0 \end{bmatrix}.
\end{align*}
(Here $\nu$ is the Lagrange multiplier associated with the initial condition $x_0 = x_{\rm init}$.)
\fi

\section{Battery control example}
\label{s-examples}

We now demonstrate the RS-MPC policy 
on an example of controlling a battery to power an uncertain load
while minimizing the cost of grid power.

\ifieee
\else
    (See figure~\ref{f-microgrid}.)

    \begin{figure}
    \centering
    \ifplots
    \begin{tikzpicture}

    \tikzset{
      position/.style args={#1:#2 from #3}{
        at=(#3.#1), anchor=#1+180, shift=(#1:#2)
      }
    }

    \tikzstyle{device}=[draw, semithick, fill=gray,
                        node distance=3cm, minimum height=0.7cm, minimum width=1.5cm]
    \tikzstyle{net}=[draw, circle, semithick, node distance=3cm,  
                     minimum width=.4cm, inner sep=0cm]
    \tikzstyle{classicterm} = [draw, semithick]
    \tikzstyle{terminal} = [draw, semithick, -latex']
    \tikzstyle{termnoarrow} = [draw, semithick]
    \tikzstyle{bus} = [draw, line width=3pt]

    \node (net)  [net] {};
    \node (battery)  [device, above of=net, node distance=1.8cm] {battery};
    \node (grid) [device, draw=none, minimum height=0.7cm, fill=none, left of=net, node distance=3.5cm] {grid};
    \node (load) [device, draw=none, minimum height=0.7cm, fill=none, right of=net, node distance=3.5cm] {load};

    \node (rightend) [position=000:{3cm} from net] {};
    \node (leftend)  [position=180:{3cm} from net] {};

    \draw [semithick, -latex'] ([xshift=3pt] battery.south) -- ([xshift=3pt] net.north)
             node[midway, right] {$p^{\rm batt}_t$};
    \draw [semithick, -latex'] ([xshift=-3pt] net.north) -- ([xshift=-3pt] battery.south);
    \draw [semithick, -latex'] (grid) -- (net) node[midway, below] {$p^{\rm grid}_t$};
    \draw [semithick, -latex'] (net) -- (load) node[midway, below] {$p^{\rm load}_t$};
    \draw [semithick] (grid.east) -- (grid.15);
    \draw [semithick] (grid.east) -- (grid.345);
    \draw [semithick] (load.west) -- (load.165);
    \draw [semithick] (load.west) -- (load.195);

    \end{tikzpicture}
    \fi
    \caption{Battery charge control schematic.}
    \label{f-microgrid}
    \end{figure}
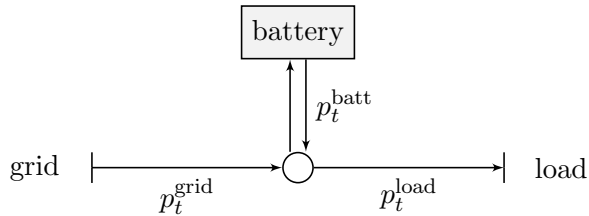
\fi

\subsection{Model}

\paragraph{Battery\pp}
In time period $t$, the battery charge 
is $q_t$ and the discharge power is $p^{\rm batt}_t$.
The battery dynamics are
\[
q_{t+1} = q_t - hp^{\rm batt}_t 
 \quad t=0, \dots, T,
\]
where $h$ is the length of a single time interval.
The battery charge must satisfy $0 \le q_t \le q^{\rm max}$
and the initial condition $q_0 = q_{\rm init}$.
Here $p_t^{\rm batt}$ is the amount of power discharged from the battery at time $t$.

\paragraph{Grid connection\pp}
The power from the grid at time $t$ is $p^{\rm grid}_t$.
For each unit of energy purchased from the grid at time $t$,
we pay $c_t$ dollars;
the total cost is $h \sum_{t=0}^{T-1} c_t p^{\rm grid}_t$.

\paragraph{Net power demand\pp}
The load demand at time $t$ is $p^{\rm load}_t$.
We assume the load is net of any solar or wind generation,
and may therefore be negative.
The load power demand must be met at every time period, \ie,
\[
p^{\rm load}_t \le p^{\rm grid} + p^{\rm batt}.
\]
The net power demand is a stochastic process described by the first-order auto-regressive model
\begin{equation}
\label{e-pload-dynamics}
p^{\rm load}_{t+1} = \alpha p^{\rm load}_t + (1-\alpha) p^{\rm base}_t + w_t.
\end{equation}
Here $p^{\rm base}_t$ is the baseline power demand,
\ie, it is the typical demand that would be expected at time $t$
in the absence of additional information.
The coefficient $\alpha > 0$ models reversion of the demand power to the baseline value.
The noise $w_t \sim \mathcal N(0, \sigma^2)$ is Gaussian and independent across time periods,
with rate function is $\rho(w) = w^T w / (2\sigma^2)$
and cumulant generating function $c(z) = (\sigma^2/2) z^T z$.
This type of auto-regressive model with a baseline value is common;
see \cite[\S A]{moehle2019dynamic} for details.

\paragraph{Prescient problem\pp}
The problem of minimizing the cost of grid power
can be cast as a linear convex stochastic control problem%
\ifieee
    , but we do not give it here.
\else
    \ifacc
        , but we do not give it here.
    \else
        ; the exact parameterization is given in appendix~\ref{s-battery-parametrization}.
    \fi
\fi
Here we simply note that the prescient problem \eqref{e-prec-prob} is
\begin{equation} \label{e-prec-prob-battery}
\begin{array}{ll}
\mbox{minimize} 
    & h \sum_{t=1}^{T-1} c_t p^{\rm grid}_t \\
\mbox{subject to}  
    & q_{t+1} = q_t - h p^{\rm batt}_t%
        \ifdoublecolumn \\ \else , \quad t=0, \dots, T-1 \\ \fi
    \ifdoublecolumn
        & 0 \le q_{t+1}\le q^{\rm max} \\
    \else 
        & 0 \le q_t \le q^{\rm max}, \quad t=1, \dots, T \\ 
    \fi
    & q_0 = q_{\rm init} \\
    & p^{\rm load}_{t+1} = \alpha p^{\rm load}_t + (1-\alpha) p^{\rm base}_t + w_t
        \ifdoublecolumn \\ \else , \quad t=0, \dots, T-1 \\ \fi
    & p^{\rm load}_t \le p^{\rm batt}_t + p^{\rm grid}_t
        \ifdoublecolumn \\ \else , \quad t=0, \dots, T-1 \\ \fi
    & p^{\rm grid}_t \ge 0.
\end{array}
\end{equation}
The variables are $q_t$, for $t = 0, \dots, T$,
as well as $p^{\rm grid}_t$, $p^{\rm batt}_t$, and $p^{\rm load}$, for $t = 0, \dots, T-1$.
\ifdoublecolumn
    All constraints indexed by $t$ hold for $t = 0,\dots, T-1$.
\else
\fi

\paragraph{Algorithm interpretation\pp}
To carry out one iteration in the convex--concave procedure,
we first solve the prescient problem~\eqref{e-prec-prob-battery},
then set $w^{(k)}$ to be the gradient of the cumulant generating function 
at the optimal dual variables $\lambda$ of the load dynamics~\eqref{e-pload-dynamics}.
The optimal dual variable $\lambda_t$ can be interpreted as the price of energy at time $t$
\cite[\S 2.3]{moehle2019dynamic}.
This means that RS-MPC
pessimistically assumes there will be greater demand precisely when the price of energy is high.


\paragraph{Data\pp}
We used the parameter values
$q^{\rm init} = 2.5$ kWh,
$q^{\rm max} = 5$ kWh,
and
$\alpha = 0.5$.
The planning horizon is $T = 300$ time steps,
with the discretization interval $h$ chosen so that the planning horizon $hT$ is two days.
The price of energy $c$ is
\[
c_t = \begin{cases}
\text{$15$ \cent/kWh} & \text{$t$ is between 21:00 and  6:00} \\
\text{$40$ \cent/kWh} & \text{$t$ is between 13:00 and 19:00} \\
\text{$25$ \cent/kWh} & \text{otherwise.} \\
\end{cases}
\]
The baseline load $p^{\rm bl}_t$ is shown 
along with the results in figure~\ref{f-battery-trajectories}.
Note that power demand is low in the morning, negative in the afternoon (due to solar generation),
and high in the evening.

%
%

\subsection{Results}

\begin{figure}
\centering
\ifplots
\begin{tikzpicture}
    \begin{groupplot}[
        group style = {group size=1 by 3, vertical sep=0.3cm},
        table/col sep=comma,
        xmin = 0,
        xmax = 48,
        width = 1.0\columnwidth,
        height = 0.4\columnwidth,
    ]

    \nextgroupplot[
        xticklabels={,,},
        ylabel={$p_t^{\rm batt}$, $p_t^{\rm load}$ (kW)},
    ]
        \addplot [black, semithick        ] table[x=times, y= P_req_rn]{\datadir/battery_trajectories.\ext};
        \addplot [black, semithick, dashed] table[x=times, y= P_req_ra]{\datadir/battery_trajectories.\ext};
        \addplot [ blue, semithick        ] table[x=times, y=P_grid_rn]{\datadir/battery_trajectories.\ext};
        \addplot [green, semithick        ] table[x=times, y=P_grid_ra]{\datadir/battery_trajectories.\ext};
        \addplot [  red, semithick        ] table[x=times, y=P_grid_cl]{\datadir/battery_trajectories.\ext};

    \nextgroupplot[
        xticklabels={,,},
        ylabel={$q_t$ (kWh)},
    ]
        \addplot [ blue, semithick] table[x=times, y=E_rn]{\datadir/battery_trajectories.\ext};
        \addplot [green, semithick] table[x=times, y=E_ra]{\datadir/battery_trajectories.\ext};
        \addplot [  red, semithick] table[x=times, y=E_cl]{\datadir/battery_trajectories.\ext};

    \nextgroupplot[
        ylabel={$\lambda_t$ (\$/kWh)},
        xlabel = Time (h),
    ]
        \addplot [black, semithick] table[x=times, y=   price]{\datadir/battery_trajectories.\ext};
        \addplot [ blue, semithick] table[x=times, y=price_rn]{\datadir/battery_trajectories.\ext};
        \addplot [green, semithick] table[x=times, y=price_ra]{\datadir/battery_trajectories.\ext};
        \addplot [  red, semithick] table[x=times, y=price_cl]{\datadir/battery_trajectories.\ext};

    \end{groupplot}
\end{tikzpicture}
\fi
\caption{
Three trajectories for the battery control example:
risk-neutral control under the most likely outcome $w = \Expect w = 0$ (blue),
risk-averse control with the projected unfavorable outcome for $w$ (green),
and risk-averse control under the most likely outcome $w = \Expect w = 0$ (red).
The solid black curve shows the mean load power (with $w = 0$),
and the dashed black curve shows the power demand under the unfavorable outcome.
}
\label{f-battery-trajectories}
\end{figure}
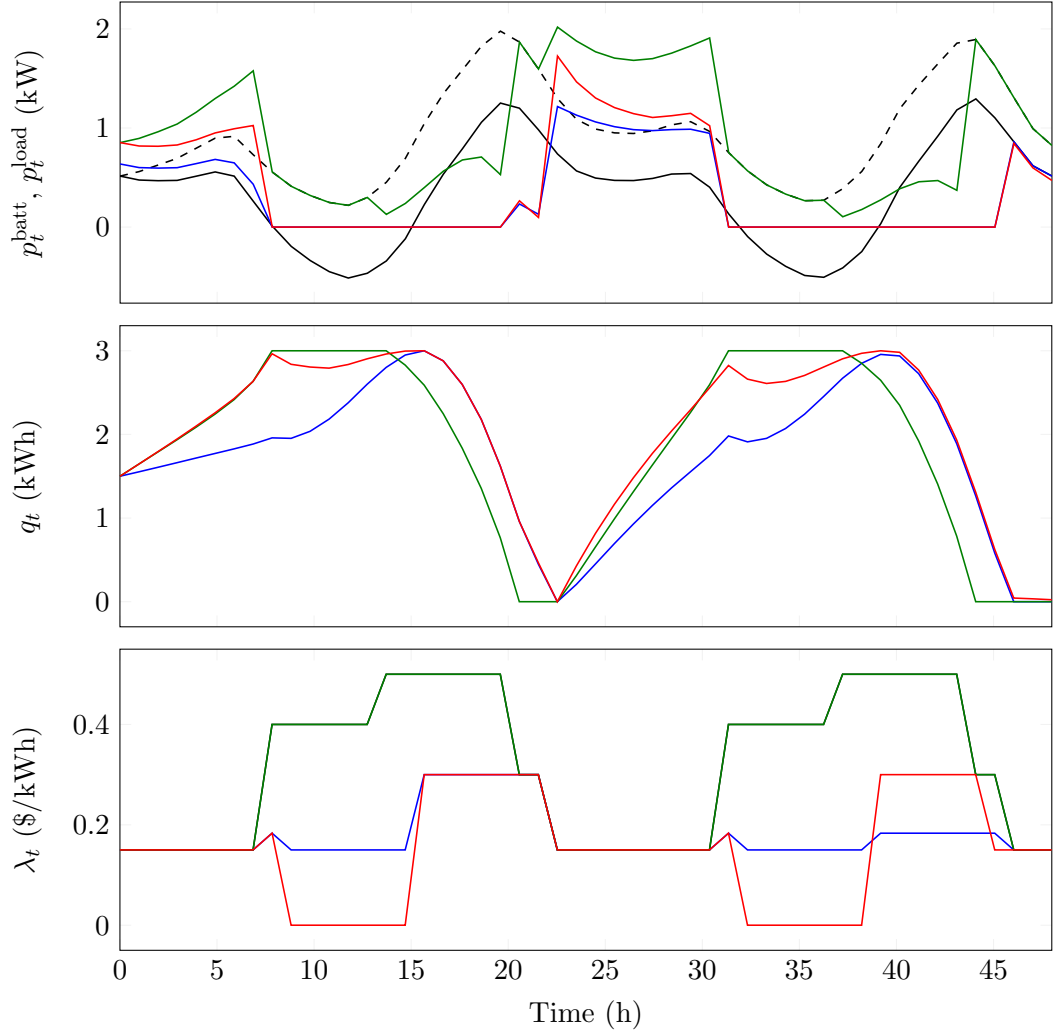

\paragraph{Trajectory comparison\pp}
In figure~\ref{f-battery-trajectories},
we show three sets of trajectories for the battery charge control problem.
Each set consists of the grid power consumption (top plot),
the battery charge (middle plot),
and the price of energy 
\ie, the optimal dual variable for constraint~\eqref{e-pload-dynamics}
(bottom plot).

In blue, we plot the optimal trajectories for the prescient problem~\eqref{e-prec-prob}
with realized outcome $w = \Expect w = 0$.
(This trajectory would be used by risk-neutral MPC
to choose the first control input.)
This plan begins charging the battery in the morning,
relying on afternoon solar power to finish charging.
The battery is discharged in the evening
when grid power is expensive and the demand is high.
The local price of energy is flatter than the grid price,
because we use the battery to shift our power purchases to be earlier in the day.

In green, we plot the optimal trajectory for \eqref{e-prec-prob},
where $w$ is chosen adversarially,
\ie, it optimizes the bound~\eqref{e-bound-risk-averse} with $\gamma = 2$.
(This trajectory would be used by RS-MPC.)
This plan charges the battery completely in the morning,
because it assumes no excess solar production in the afternoon.
The local price of energy is higher than in the case $w = 0$,
because we pessimistically assume higher power demand, especially during peak hours.

Finally, the trajectory in red is a closed-loop simulation of RS-MPC
under the outcome $w = \Expect w = 0$.
This means that although the policy is planning against an adversarial outcome,
the true outcome is not chosen adversarially.
This allows us to compare RS-MPC against the optimal prescient plan
for this particular outcome (shown in blue).
Because of our pessimism,
we charge more aggressively in the morning than the blue (risk-neutral) trajectory,
because we are planning for higher demand throughout the day.
Because the true outcome is $w = \Expect w = 0$,
this pessimism is misplaced (in this particular example),
and the local price of energy is more uneven than for the optimal risk-neutral trajectory,
\ie, RS-MPC produces more price fluctuations.
This is because the policy has saved too much energy in the morning,
and has a surplus in the afternoon, causing the price to decrease.

\begin{figure}
\ifplots
\centering
\begin{tikzpicture}
    
    \newcommand\xmin{0}
    \newcommand\xmax{7}

    \begin{groupplot}[
        group style = {group size=1 by 3, vertical sep=0.3cm},
        table/col sep=comma,
        xmin = \xmin,
        xmax = \xmax,
        ymin = 0,
        ymax = 550,
        width = 1\columnwidth,
        height = 0.30\columnwidth,
        yticklabels={,,},
    ]

    \nextgroupplot[
        xticklabels={,,},
    ]
        \addplot [
            hist={data=x, bins=50, data min=\xmin, data max=\xmax},
            fill=gray,
        ] table[x index=0] {\datadir/battery_obj_vals.\ext};
        \addplot [     blue, semithick] table[y index=0, x index=1]{\datadir/battery_vert_line_data.\ext};
        \addplot [green, semithick] table[y index=0, x index=2]{\datadir/battery_vert_line_data.\ext};
        \addplot [      red, semithick] table[y index=0, x index=3]{\datadir/battery_vert_line_data.\ext};

    \nextgroupplot[
        xticklabels={,,},
    ]
        \addplot [
            hist={data=x, bins=50, data min=\xmin, data max=\xmax},
            fill=gray,
        ] table[x index=1] {\datadir/battery_obj_vals.\ext};
        \addplot [     blue, semithick] table[y index=0, x index=4]{\datadir/battery_vert_line_data.\ext};
        \addplot [green, semithick] table[y index=0, x index=5]{\datadir/battery_vert_line_data.\ext};
        \addplot [      red, semithick] table[y index=0, x index=6]{\datadir/battery_vert_line_data.\ext};

    \nextgroupplot[
        xlabel = Closed-loop cost $C_\pi(w)$,
    ]
        \addplot [
            hist={data=x, bins=50, data min=\xmin, data max=\xmax},
            fill=gray,
        ] table[x index=2] {\datadir/battery_obj_vals.\ext};
        \addplot [     blue, semithick] table[y index=0, x index=7]{\datadir/battery_vert_line_data.\ext};
        \addplot [    green, semithick] table[y index=0, x index=8]{\datadir/battery_vert_line_data.\ext};
        \addplot [      red, semithick] table[y index=0, x index=9]{\datadir/battery_vert_line_data.\ext};

    \end{groupplot}
\end{tikzpicture}
\fi
\caption{
The distribution of costs $C_\pi$ obtained using the RS-MPC
policy with risk aversion parameter $\gamma = 0$ (top),
$\gamma = 2$ (middle), and $\gamma = 5$ (bottom).
The vertical lines show the values of $J_\pi = R_\gamma(C_\pi)$ obtained in closed loop,
evaluated for all three values of $\gamma$
($\gamma = 0$ in blue, $\gamma = 2$ in green, and $\gamma = 5$ in red).
}
\label{f-battery-histograms}
\end{figure}
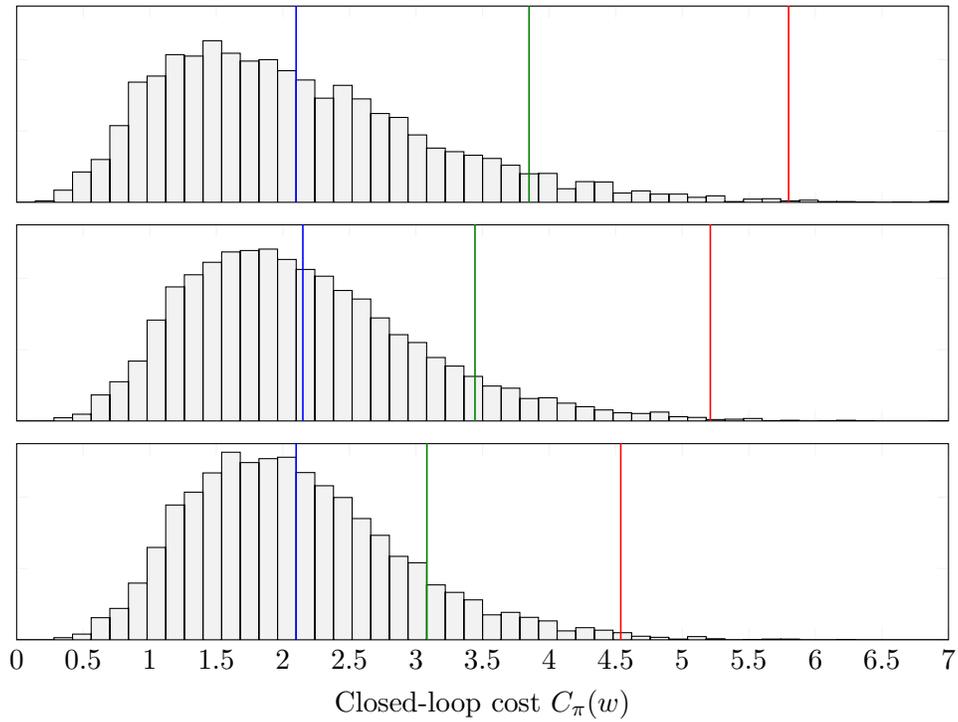

\paragraph{Cost distribution\pp}
In figure~\ref{f-battery-histograms},
we show the distribution of costs $C_\pi(w)$
achieved for risk-neutral MPC ($\gamma = 0$)
and RS-MPC ($\gamma = 2$ and $\gamma = 5$).
We observe that RS-MPC reduces the probability of achieving a very high cost.

We also show the risk-adjusted cost $J_\pi = R_\gamma(C_\pi(w))$,
obtained in closed loop, for all three values of $\gamma$.
RS-MPC reduces $J$ when $\gamma$ is high,
\ie, when the `true' cost is risk averse.
When the true cost is risk-neutral, \ie, when $J$ is evaluated using $\gamma = 0$
(shown by vertical blue lines in figure~\ref{f-battery-histograms}),
we observe a surprising result:
the performance of RS-MPC is comparable to risk-neutral MPC.
We suspect the cautious planning of RS-MPC
helps avoid being caught mid-day with little battery charge,
and therefore having to purchase grid power when it is most expensive.
(This phenomenon does not hold for all examples;
for example, for the LQR problem, risk-neutral MPC is in fact optimal,
and risk-averse policies are typically suboptimal when the true cost is risk neutral.)

\ifieee
\else
\section{Conclusion}
In this paper, we address risk-sensitive convex stochastic control problems
by approximating them as deterministic optimization problems
In future work, we hope to expand the set of problems that can be addressed by these techniques.
In particular, we will
minimize a sum of exponentials of convex stage costs
instead of an exponential of a sum of convex stage costs.
This allows us to consider other interesting risk-averse problems,
such as the Merton's consumption--investment problem.
\fi

\paragraph{Acknowledgments\pp}
I would like to thank Stephen Boyd for useful discussions and feedback.


\ifelsevier
    \bibliographystyle{plainnat}
    \bibliography{risk_averse_ctrl}
\else
    \printbibliography
    \ifieee
    \else
        \ifacc
        \else
            \clearpage
        \fi
    \fi
\fi

\appendix

\ifieee
\else
    \section{Proof of risk bound}
\fi

\label{s-varadhan-bound}

Here we prove inequality~\eqref{e-varadhan-bound}.


\paragraph{Affine functions\pp}
First consider the case that $f$ is affine, \ie, $f(z) = a^T z + b$.
From the definitions of the risk operator \eqref{e-risk-operator}
and cumulant generating function \eqref{e-cgf},
it can be verified that
\[
R_\gamma \big(f(w)\big) = \frac1\gamma c(\gamma a) + b.
\]
Because the cumulant generating function is the conjugate of the rate function,
this is
\begin{equation}
\label{e-varadhan-affine}
\begin{aligned}
\ifdoublecolumn
    R_\gamma \big(f(w)\big) 
    &= \frac1\gamma \sup_z \big( \gamma a^T z + \gamma b - \rho(z)\big) \\
    &= \frac1\gamma \sup_z \big( \gamma f(z) - \rho(z)\big),
\else
    R_\gamma \big(f(w)\big) 
    = \frac1\gamma \sup_z \big( \gamma a^T z + \gamma b - \rho(z)\big) 
    = \frac1\gamma \sup_z \big( \gamma f(z) - \rho(z)\big),
\fi
\end{aligned}
\end{equation}
\ie, the bound holds with equality for affine functions.

\paragraph{Convex functions\pp}
If $f$ is convex,
we apply \eqref{e-varadhan-affine} to $\hat f$, a first-order Taylor expansion of $f$
around a maximizing value of $z$,
such that $\hat f \le f$ and 
\begin{align}
\label{e-varadhan-functions-are-equal}
\sup_z \big( \gamma \hat f(z) - \rho(z)\big)
=
\sup_z \big( \gamma f(z) - \rho(z)\big).
\end{align}
(If no such maximizer $z$ exists, $\hat f$ is a limit of Taylor expansions around a sequence
of points that attain the supremum in the limit.)
From this we obtain
\begin{align*}
R_\gamma \big(f(x)\big) 
&\ge R_\gamma \big(\hat f(x)\big) \\
&= \frac1\gamma \sup_z \big( \gamma \hat f(z) - \rho(z)\big) \\
&= \frac1\gamma \sup_z \big( \gamma f(z) - \rho(z)\big).
\end{align*}
The first line follows from $\hat f \le f$ and the apparent monotonicity of the risk operator,
the second line from \eqref{e-varadhan-affine} applied to the affine function $\hat f$,
and the third line from \eqref{e-varadhan-functions-are-equal}.

\ifieee
\else
    \ifacc
    \else
        \section{Parameterization of battery example}
        \label{s-battery-parametrization}
        We can express the battery charge control problem as a linear convex stochastic control problem with dynamics
        given by \eqref{e-dynamics} with
        state $x_t = (q_t, p^{\rm load}_t)$,
        input $u_t = (p^{\rm batt}_t, p^{\rm grid}_t)$,
        and noise $w_t' = (1-\alpha) p^{\rm base}_t + w_t$.
        The dynamics parameters are
        \[
        A_t = \begin{bmatrix} 1 & 0 \\ 0 & \alpha \end{bmatrix}, \quad
        B_t = \begin{bmatrix} -h & 0 \\ 0 & 0 \end{bmatrix},
        \]
        and the stage cost functions are

        \[
        g_t(x_t, u_t) 
        = \begin{cases}
        hc_t p^{\rm grid}_t & \ifdoublecolumn
                                  \begin{aligned}
                                      & p^{\rm load}_t \le p^{\rm grid}_t + p^{\rm batt}_t,  \\
                                      & 0 \le q \le q^{\rm max},\; p^{\rm grid}_t \ge 0
                                  \end{aligned}
                                  \\
                              \else
                                  p^{\rm load}_t \le p^{\rm grid}_t + p^{\rm batt}_t,\;
                                  0 \le q \le q^{\rm max},\;
                                  p^{\rm grid}_t \ge 0 \\
                              \fi
        \infty & \text{otherwise.}
        \end{cases}
        \]
    \fi
\fi

%
%
%
%
%
%
%
%
%
%
%
%
%
%
\end{document}